\documentclass[12pt]{article}
\usepackage[top=25.4mm, bottom=25.4mm, left=25.4mm, right=25.4mm]{geometry}

\usepackage{amssymb,amstext,amsmath}
\usepackage{graphicx}% Include figure files
\usepackage{dcolumn}% Align table columns on decimal point
\usepackage{bm} % bold math
\usepackage{multirow}
\usepackage{bigstrut}
\usepackage{makecell,rotating}
\usepackage{mathrsfs}
\usepackage{booktabs}
\usepackage{threeparttable}
\usepackage{chngpage}
\usepackage{float}
\usepackage{color}
\usepackage{caption,subcaption}
\usepackage{scicite}
\usepackage{xr}
\newcommand{\x}{\mathbf{x}}
\newcommand{\U}{\mathbf{u}}
\newcommand{\A}{\mathbf{A}}
\newcommand{\B}{\mathbf{B}}
\newcommand{\G}{\mathbf{G}}
\newcommand{\p}{\mathbf{P}}
\newcommand{\Q}{\mathbf{Q}}
\newcommand{\M}{\mathbf{M}}
\newcommand{\I}{\mathbf{I}}
\newcommand{\F}{\mathbf{F}}
\newcommand{\X}{\mathbf{\Xi}}

\newcommand{\ie}{\textit{i.e.}}

\usepackage{color}
	 \definecolor{darkred}{rgb}{0.75,0,0}
	 \definecolor{darkgreen}{rgb}{0,0.5,0}
	 \definecolor{darkblue}{rgb}{0,0,0.75}
  	 \definecolor{darkorange}{rgb}{1,0.9,0.1}
	 \definecolor{dark}{rgb}{0,0,0}

\begin{document}

\title{Energy cost for controlling complex networks}
%\title{Minimum energy control in complex networks}
\author{Gaopeng Duan$^{1,2^*}$, Aming Li$^{3,4^*}$, Tao Meng$^{1}$, Guofeng Zhang$^{2}$ and Long Wang$^{1^\dag}$}
\date{\today}
\maketitle
\begin{enumerate}
  \item Center for Systems and Control, College of Engineering, Peking University, Beijing 100871, China
    \item Department of Applied Mathematics, Hong Kong Polytechnic University, Hong Kong
   \item Department of Zoology and Oxford Centre for Integrative Systems Biology, University of Oxford, Oxford OX1 3PS, UK
      \item Chair of Systems Design, Department of Management, Technology and Economics, ETH Z\"urich, Weinbergstrasse 56/58, Z\"urich CH-8092, Switzerland
   \item[$*$]  These authors contributed equally to this work
   \item[$\dag$] Correspondence to: longwang@pku.edu.cn
\end{enumerate}

\begin{abstract}
The controllability of complex networks has received much attention recently, which tells whether we can steer a system from an initial state to any final state within finite time with admissible external inputs.
In order to accomplish the control in practice at the minimum cost, we must study how much control energy is needed to reach the desired final state.
At a given control distance between the initial and final states, existing results present the scaling behavior of lower bounds of the minimum energy in terms of the control time analytically.
However, to reach an arbitrary final state at a given control distance, the minimum energy is actually dominated by the upper bound, whose analytic expression still remains elusive.
Here we theoretically show the scaling behavior of the upper bound of the minimum energy in terms of the time required to achieve control.
Apart from validating the analytical results with numerical simulations, our findings are feasible to the scenario with any number of nodes that receive inputs directly and any types of networks.
Moreover, more precise analytical results for the lower bound of the minimum energy are derived in the proposed framework.
Our results pave the way to implement realistic control over various complex networks with the minimum control cost.
\end{abstract}

%\vspace{2cm}
%\newpage
%\tableofcontents

\section{Introduction}

An ultimate goal of studying complex systems is to control them on the basis of the underlying topological structures, where nodes indicate units of a system and edges capture who interacts with whom \cite{Liu2016Rev,barabasi2016network,Liu2008Controllability,Xie2002,Xie2003}.
Indeed, by implementing appropriate external control signals, if we can drive a system from an arbitrary initial state to any final state in finite time, we say that the system is controllable, \ie, in principle, we are able to steer the system along our expectations.
Recently, the problem of finding set of minimal number of nodes that receive external inputs directly to make a network controllable has been investigated\cite{Liu2011,Wang2013}.
And in the past several years, several important results have elucidated important problems pertaining to
%\cite{Pequito2013,Olshevsky2014},
node classification\cite{Jia2013,Vinayagam03052016}, control profiles\cite{Ruths2014science}, target control\cite{Gao2014}, control of edge dynamics\cite{Nepusz2012}, as well as the energy (or cost) required for control\cite{Yan2012PRL,energy2014,Yan2015a,Chen2016Energy,Li2017,Li2017ConEng}.

Beyond the basic property, namely controllability of a system, the control energy steering the system from an initial to a final state has received much attention recently. Indeed, the energy
%evaluates not only the strength of external inputs but also the time to reach the final state,
tells the cost required to pay in practical control, and thus represents another dimension of difficulty in achieving control.
Although theoretically approximate lower bound of control energy and its scaling behavior in terms of the control time have been provided in the literarure for both static and temporal networks, the energy to reach an arbitrary final state in phase space is usually dominated by the upper bound \cite{Yan2012PRL,Li2017ConEng}.
Analytical forms on the upper bound of control energy are as yet still missing, and the existing results are all based on the myriad numerical calculations.
%And the lower bound of energy is also employed to explore how the brain mores between cognitive states \cite{Gu08Con,Bassett18NP}.
%But they are mere for extreme values of the energy among all final states which are equal to the initial state, namely they are the minimum energy corresponding to two specific final states.
%Here for finite control time, we systematically investigate the minimum control energy for specific final state.
In this article, apart from presenting more precise lower bound of the minimum control energy, we theoretically derive the upper bound for the first time.
Furthermore, we show the scaling behavior of both bounds, and numerical validations are also given for both cases.
\section{The minimum energy for controlling complex \text{networks}}
Here we consider the canonical linear time-invariant dynamics
\begin{equation}\label{sys1}
\dot{\x}(t)=\A\x(t)+\B\U(t),
\end{equation}
where $\x(t)=(x_1(t)~ x_2(t)~ \dots ~x_n(t))^\text{T}$ is the state of the whole network with $x_i(t)$ capturing the state of node $i$; $\U(t)=(u_1(t)~ u_2(t)~ \dots ~u_m(t))^\text{T}$ is the control input; $\A=(a_{ij})_{n n}$ is the adjacent matrix of the network; $\B=(b_{ij})_{n m}$ is the input matrix with size $n\times m$, and the entry at row $i$ and column $j$ is $b_{ij}$, being $1$ if node $i$ receives the external control input signal $u_j(t)$ directly (driver node), being 0 otherwise.

The networked system (\ref{sys1}) is said to be controllable, if it can be driven from any initial state $\x_0=\x(t_0)$ toward any target state $\x_f=\x(t_f)$ at a given control time $t_f$, and the corresponding input control energy cost is defined as $E(t_0, t_f)=\int^{t_f}_{t_0}\|\U(t)\|^2\text{d}t$ with $\|\U(t)\|$ being the Euclidean norm of the vector $\U(t)$.
To minimize the above energy cost, one can adopt the minimum energy control input $\U^*(t)=\B^\text{T}\text{e}^{\A^\text{T}(t_f-t)}\G^{-1}\delta$ with $\G=\int^{t_f}_{t_0}\text{e}^{\A (t-t_0)}\B\B^\text{T}\text{e}^{\A^\text{T}(t-t_0)}\text{d}t$ and $\delta=\x_{f}-\text{e}^{\A t_f}\x_0$ \cite{OptimalBooLewis}, which gives the minimum energy cost $E(t_f)=\delta^\text{T}\G^{-1}\delta$ from $\x_0$ to $\x_f$. By assuming $t_0=0$ and $\x_0=\mathbf{0}$ for simplicity, we obtain  the minimum energy
\begin{equation}\label{E}
E(t_f)=\x^\text{T}_f\G^{-1}\x_f,
\end{equation}
and note that here the matrix $\G$ is positive definite when system (\ref{sys1}) is controllable \cite{Kalman63}.
Note that when we refer to control energy later, we mean the minimum control energy.
Clearly, for the normalized control distance $\|\x_f\|=1$ we have
\begin{equation}\label{EB}
\frac{1}{\lambda_{\max}(\G)}\leq E(t_f)\leq\frac{1}{\lambda_{\min}(\G)}.
\end{equation}
In what follows, for ease of presenting our framework, we consider undirected networks, where $\A$ corresponds to the real symmetric matrix. Subsequently, we have $\A=\p\X\p^\text{T}$ with $\p\p^\text{T}=\p^\text{T}\p=\mathbf{I}$, where $\X=\text{diag}(\lambda_1, \lambda_2, \dots, \lambda_n)$, and $\lambda_i, (i=1, 2, \dots, n)$ is the eigenvalue of $\A$ with the ascending order $\lambda_1\leq\lambda_2\leq\dots\leq\lambda_n$.
By letting $\Q=\p^\text{T}\B\B^\text{T}\p=(q_{ij})_{n n}$ and $\F=(f_{ij})_{n n}$ with $f_{ij}=\frac{1}{\lambda_i+\lambda_j}\left[\text{e}^{(\lambda_i+\lambda_j)t_f}-1\right]$, we have $\int^{t_f}_0\text{e}^{\X t}\p^\text{T}\B\B^\text{T}\p \text{e}^{\X t}\text{d}t=(q_{ij} f_{ij})_{n n}$.
Note that the limit of $f_{ij}$ is $t_f$ as $\lambda_i+\lambda_j\rightarrow0$, which keeps the above expression of $f_{ij}$ alive when $\lambda_i+\lambda_j=0$.
%if $\lambda_i+\lambda_j\neq 0$, then
%\begin{align}
%\int^{t_f}_0\text{e}^{(\lambda_i+\lambda_j)t}\text{d}t&=\frac{1}{\lambda_i+\lambda_j}\text{e}^{(\lambda_i+\lambda_j)t}\left|^{t_f}_0\right.\notag\\
%&=\frac{1}{\lambda_i+\lambda_j}\left(\text{e}^{(\lambda_i+\lambda_j)t_f}-1\right)\label{mij1};
%\end{align}
%if $\lambda_i+\lambda_j=0$, then
%\begin{equation}\label{mij2}
%\int^{t_f}_0\text{e}^{(\lambda_i+\lambda_j)t}\text{d}t=t_f.
%\end{equation}
%Since $\lim_{\lambda_i+\lambda_j\rightarrow0}\frac{\text{e}^{(\lambda_i+\lambda_j)t_f}-1}{\lambda_i+\lambda_j}=t_f$, in the following calculation and analysis, we can represent $f_{ij}$ by (\ref{mij1}) in form uniformly.
%Let $\M=[m_{ij}]\in R^{n\times n}$ with $m_{ij}=q_{ij}\cdot f_{ij}$. Then we have $\G=\p\M\p^\text{T}$. Inequations (\ref{EB}) are equivalent to
%\begin{equation}\label{EB1}
%\frac{1}{\lambda_{\max}(\M)}\leq E(t_f)\leq\frac{1}{\lambda_{\min}(\M)}.
%\end{equation}
Furthermore, we can calculate $\G$ by
\begin{equation}\label{G}
\G=\p\int^{t_f}_0\text{e}^{\X t}\p^\text{T}\B\B^\text{T}\p \text{e}^{\X t}\text{d}t\p^\text{T}=\p\M\p^\text{T},
\end{equation}
where $\M=(m_{ij})_{n n}$ with $m_{ij}=q_{ij} f_{ij}$. Based on similarity between matrices $\G$ and $\M$, we know that they have the same eigenvalues. Therefore, by calculating the eigenvalues of $\M$ we can find the lower and upper bounds of the minimum energy $E(t_f)$ given in Eq.~(\ref{EB}).

\section{Results}
As discussed in the previous section, driver nodes are nodes who receive external control inputs directly. In this section, for different numbers of driver nodes, we derive the analytical bounds of the control energy separately.
For simplicity, here we assume that each single input only injects on a single driver node, and each node only receives an input at most.

\subsection{$n$ driver nodes}\label{nd=n}
In the case of $n$ driver nodes, \ie~all nodes receive external inputs directly, we have $m=n$, and $\B=\Q=\I$, which leads to a diagonal matrix $\M$ with $m_{ii}=f_{ii}$.
According to the magnitude of the control time $t_f$, the corresponding bounds are given as follows.

\begin{figure}
\centering
\includegraphics[width=1\textwidth]{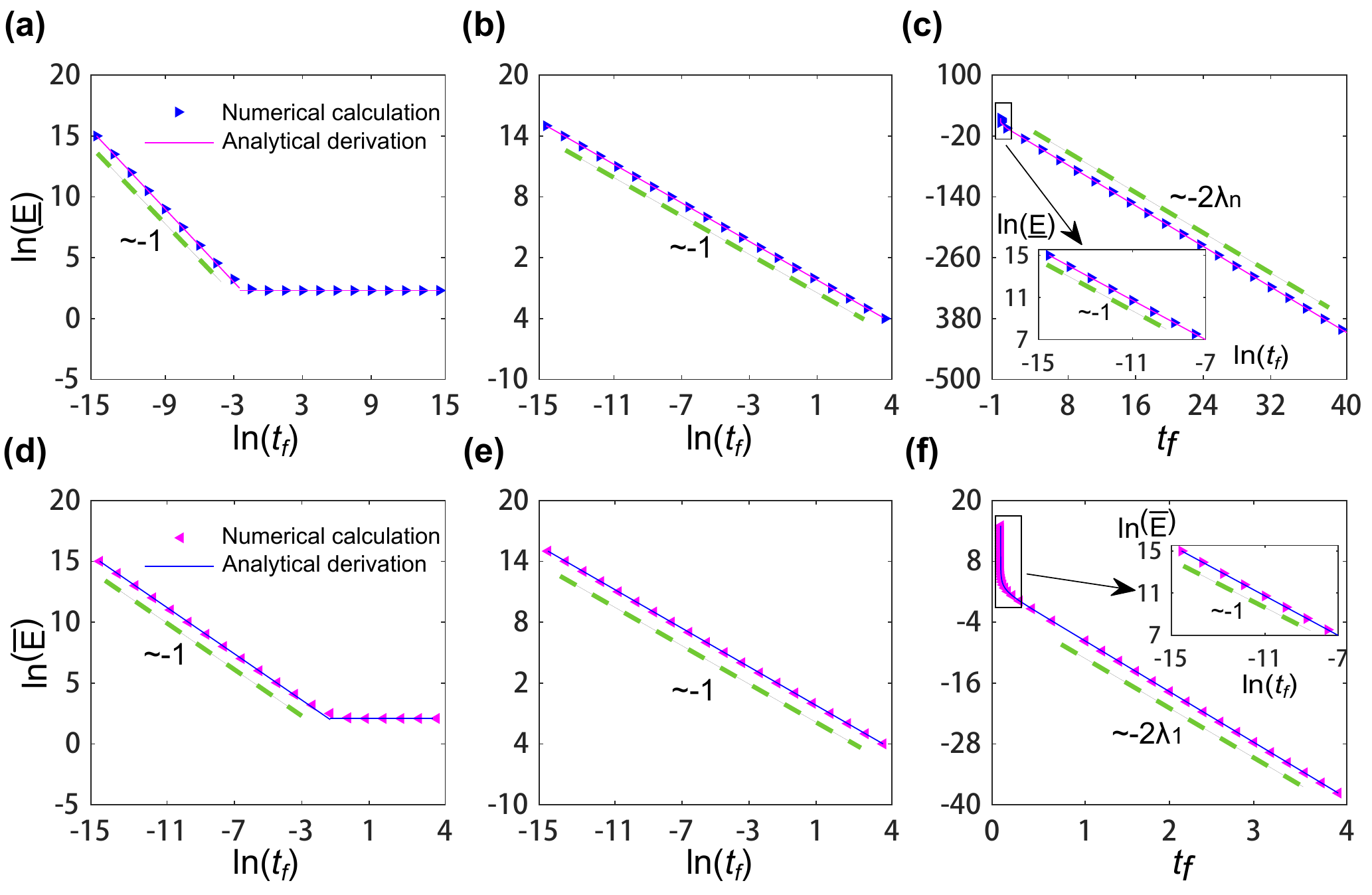}
\caption{The lower and upper bounds of control energy for $n$ driver nodes.
By controlling all nodes directly, here we show the numerical and analytical results for lower ($\underline{E}$) and upper ($\overline{E}$) bounds of control energy for different types of $\A$.
To adjust the maximum (minimum) eigenvalue of $\A$ intuitively, we set the link weight $a_{ij}$ uniformly from $[0, 1]$ in (a) to (d) and from $[-1, 0]$ in (e) and (f);
each self-loop (diagonal element) is set as $a+s_i$ with $s_i=-\sum^n_{j=1}a_{ij}$.
In (a), we set $a=-5$, which guarantees $\A$ is ND with eigenvalues in $[-14.0266, -5]$.
Similarly, in (b), $a=0$ and $\A$ is NSD with eigenvalues in $[-8.5243, 0]$.
In (c) and (d), we have $a=5$, and $\A$ is ID with eigenvalues in $[-4.0266, 5]$.
In (e), we set $a=0$, and hence $\A$ is PSD with all eigenvalues in $[0, 8.3062]$.
In (f), $a=5$ and $\A$ is PD with all eigenvalues in $[5, 13.7144]$.
In each panel, triangles (blue and purple) represent results obtained by numerical calculations and full lines indicate analytical derivations under our framework (see Sec.~\ref{nd=n} and Table~\ref{table1}).
For small $t_f$, from each panel with horizontal axis $\text{ln}(t_f)$, we see that all slopes are $-1$, which confirm our analytical results that both $\overline{E}$ and $\underline{E}$ approximate $\frac{1}{t_f}$ for different types of $\A$.
For large $t_f$, subgraphs with horizontal axis $t_f$ or ln($t_f$) show the analytical scaling behaviors of the bounds of energy precisely.
Here we adopt the BA scale-free network with $n=50$, and network is constructed based on the preferential attachment with average degree 5.8 \cite{Albert1999Diameter}.
}
\label{fig1}
\end{figure}

When $t_f$ is small, we have $\text{e}^{2\lambda_it_f}\approx1+2\lambda_it_f$, and all eigenvalues of $\M$ can be approximated by $t_f$.
Then both the upper and lower bounds of the minimum energy are $t_f^{-1}$ (see Fig.~\ref{fig1}).

When $t_f$ is large and $\A$ is indefinite (ID), \ie~$\lambda_{i-1}<0, \lambda_i=\dots=\lambda_{i+j}=0,$ $0<\lambda_{i+j+1}$, the $p$th eigenvalue of $\M$ is given by: (i) $\frac{1}{2|\lambda_p|}$ for $p=1, 2, \dots, i-1$;
(ii) $t_f$ for $p=i, i+1, \dots, i+j$;
and (iii) $\frac{\text{e}^{2\lambda_pt_f}-1}{2\lambda_p}$ for  $p=i+j+1, \dots, n$.
Therefore, we have $\lambda_{\max}(\M)=\frac{\text{e}^{2\lambda_nt_f}-1}{2\lambda_n}$ and $\lambda_{\min}(\M)\approx\frac{1}{2|\lambda_1|}$ with large $t_f$, which tells that the upper bound $\overline{E}\approx2|\lambda_1|$ and the lower bound $\underline{E}=\frac{2\lambda_n}{\text{e}^{2\lambda_nt_f}-1}\sim \text{e}^{-2\lambda_nt_f}\rightarrow0$.

Similarly,
for large $t_f$, when $\A$ is negative definite (ND, $\lambda_i<0$), $m_{ii}=\frac{\text{e}^{2\lambda_it_f}-1}{2\lambda_i}\approx\frac{-1}{2\lambda_i}$ holds. Therefore, all eigenvalues of $\M$ are approximately $\frac{1}{2|\lambda_i|}, i=1, 2, \dots, n$, respectively. Then we can obtain the upper bound of energy cost $\overline{E}\approx2|\lambda_1|$ and the lower bound of energy cost $\underline{E}\approx2|\lambda_n|$.
When $\A$ is negative semi-definite (NSD, $\lambda_{i-1}<0, \lambda_i=\dots=\lambda_n=0$), all eigenvalues of $\M$ approximate $ \frac{1}{|2\lambda_1|}, \frac{1}{|2\lambda_{2}|}, \dots, \frac{1}{|2\lambda_{i-1}|},$ $t_f, t_f, \dots, t_f$, respectively. Therefore, $\lambda_{\max}(\M)=t_f$ and $\lambda_{\min}(\M)\approx\frac{1}{2|\lambda_1|}$ with large $t_f$. Then  $\overline{E}\approx2|\lambda_1|$ and  $\underline{E}=\frac{1}{t_f}$.
When $\A$ is  positive semi-definite (PSD, $\lambda_1=\dots=\lambda_{i-1}=0, 0<\lambda_i$), all eigenvalues of $\M$ are $t_f, t_f, \dots, t_f, \frac{\text{e}^{2\lambda_it_f}-1}{2\lambda_i}, \frac{\text{e}^{2\lambda_{i+1}t_f}-1}{2\lambda_{i+1}}, \dots, $ $\frac{\text{e}^{2\lambda_nt_f}-1}{2\lambda_n}   $. Thus $\lambda_{\max}(\M)=\frac{\text{e}^{2\lambda_nt_f}-1}{2\lambda_n}\sim \text{e}^{2\lambda_nt_f}$ and $\lambda_{\min}(\M)=t_f$ for large $t_f$. Accordingly, the upper bound of energy is $\overline{E}=t^{-1}_f$ and the lower bound is $\underline{E}=\frac{2\lambda_n}{\text{e}^{2\lambda_nt_f}-1}\sim \text{e}^{-2\lambda_nt_f}$.
When $\A$ is  positive definite (PD, $0<\lambda_i$), all eigenvalues of $\M$ are $\frac{\text{e}^{2\lambda_1t_f}-1}{2\lambda_1}, \frac{\text{e}^{2\lambda_{2}t_f}-1}{2\lambda_{2}}, \dots, \frac{\text{e}^{2\lambda_nt_f}-1}{2\lambda_n}$. Obviously, $\lambda_{\max}(\M)=\frac{\text{e}^{2\lambda_nt_f}-1}{2\lambda_n}$ and $\lambda_{\min}(\M)=\frac{\text{e}^{2\lambda_1t_f}-1}{2\lambda_1}$. Consequently, $\overline{E}=\frac{2\lambda_1}{\text{e}^{2\lambda_1t_f}-1}\sim \text{e}^{-2\lambda_1t_f}$ and $\underline{E}=\frac{2\lambda_n}{\text{e}^{2\lambda_nt_f}-1}\sim \text{e}^{-2\lambda_nt_f}$.

All the above analytical scaling laws are confirmed by numerical simulations presented in Fig.~\ref{fig1}.

\subsection{One driver node}\label{nd=1}
\begin{figure}
\centering
\includegraphics[width=0.45\textwidth]{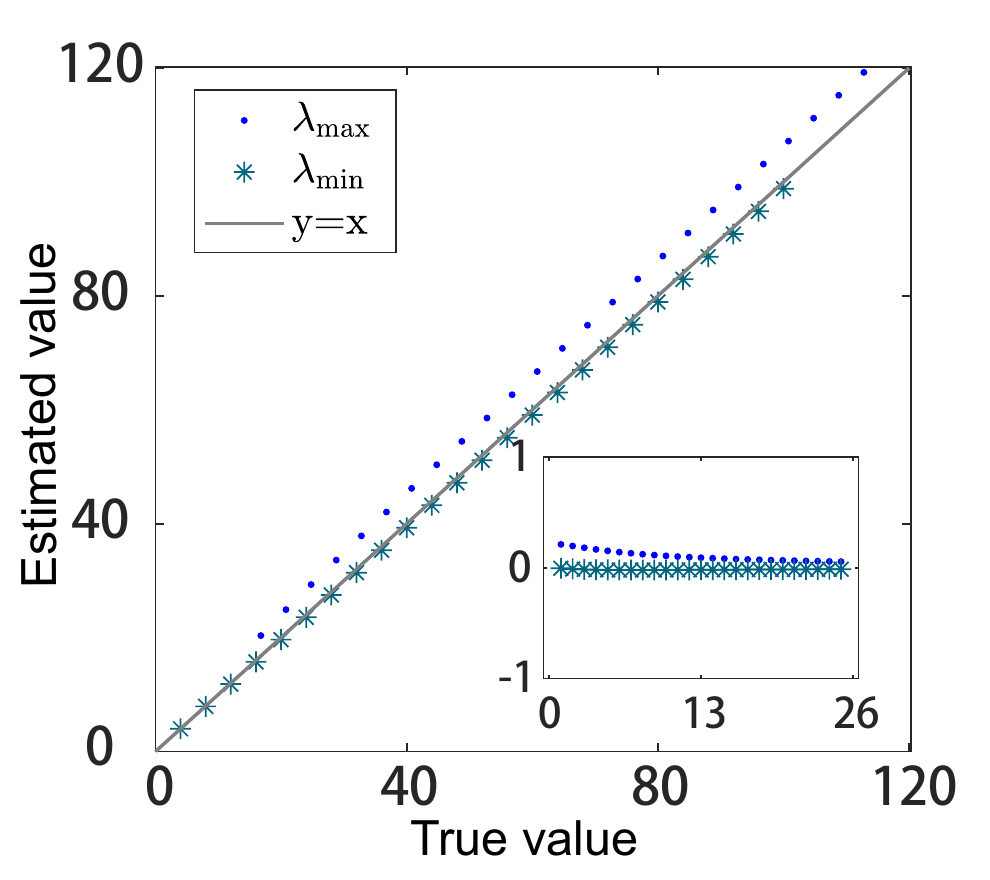}
\caption{Veracity of eigenvalues estimation based on Eqs.~(\ref{upM}) and (\ref{lowM}) for an arbitrary symmetric positive definite matrix. Here, we randomly generate $25$ matrices with minimum eigenvalue being $i\cdot 4$, $i=1, 2, \dots, 25$, where $i$ is the index of the matrix. The horizontal and vertical coordinates represent the true eigenvalues and estimated eigenvalues by Eqs.~(\ref{upM}) and (\ref{lowM}), from which it is clear the generated pattern almost overlaps with $y=x$. The inset presents ratio errors of differences between approximated eigenvalues by Eqs.~(\ref{upM}), (\ref{lowM}) and the true eigenvalues, which indicates the accuracy of estimation is reliable, especially the estimation of minimum eigenvalues by (\ref{lowM}). }\label{fig2'}
\end{figure}

In the case of one driver node, the scaling behavior of the lower bound $\underline{E}$ is given in \cite{Yan2012PRL}, in which the maximum eigenvalue of $\G$ is approximated by the trace of $\G$. In order to analytically obtain both the upper and lower bounds of the control energy $E$ shown in (\ref{EB}), we adopt the approach presented in \cite{lam2011estimates} to approximate the maximum and minimum eigenvalues of $\M$ by
\begin{equation}\label{upM}
\lambda_{\max}(\M)\approx f(\overline{\alpha}, \overline{\beta})
\end{equation}
and
\begin{equation}\label{lowM}
\lambda_{\min}(\M)\approx \frac{1}{f(\underline{\alpha}, \underline{\beta})}
\end{equation}
where $f(\alpha, \beta)= \sqrt{\frac{\alpha}{n}+\sqrt{\frac{n-1}{n}(\beta-\frac{{\alpha}^2}{n})}}$,
$
\overline{\alpha}=\text{trace}(\M^2),
\overline{\beta}=\text{trace}(\M^4),
\underline{\alpha}=\text{trace}((\M^{-1})^2),
$
and
$
\underline{\beta}=\text{trace}((\M^{-1})^4).
$
From Fig.~\ref{fig2'} we can see that it is feasible to employ (\ref{upM}) and (\ref{lowM}) to approximate respectively the maximum and the minimum eigenvalues of the real symmetric matrix with high accuracy. Specially, for positive definite matrix $\G$, the accuracy is more pronounced, as shown in Fig.~S1 in SI.

\begin{figure}
\centering
\includegraphics[width=1\textwidth]{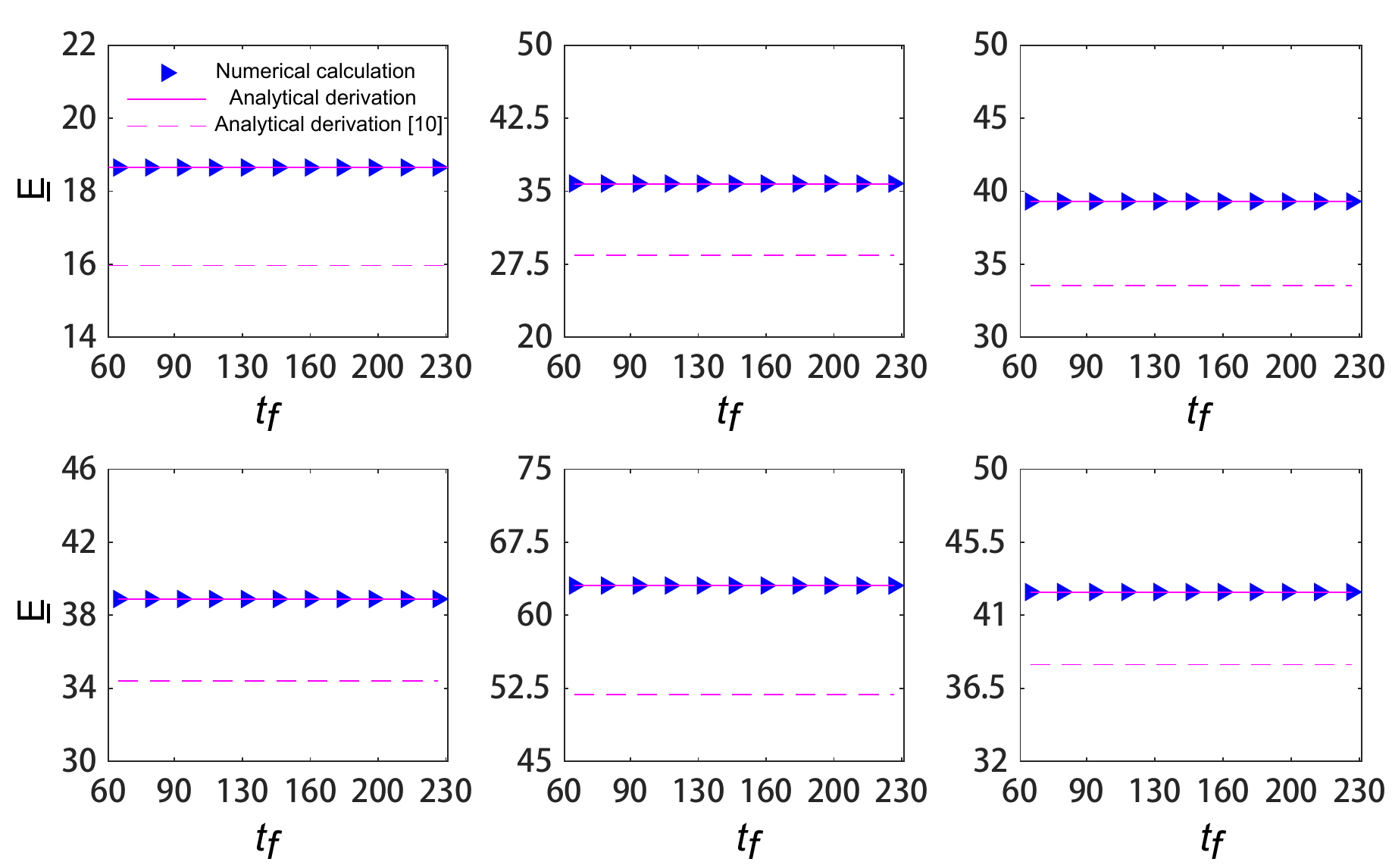}
\caption{The lower bound of energy comparisons between the methods shown in \cite{Yan2012PRL} and this paper. Here we randomly generate BA scale-free networks with $\A$ being ND (other parameters are the same as those in Fig.~\ref{fig1}) and $a_{ij}$ is selected from $[1, 3]$ uniformly with $a=-2$. For approximating the maximum eigenvalue of $\M$, here we use the method shown in (\ref{upM}), while in \cite{Yan2012PRL}, it is inferred by the corresponding trace. Since the existing results only consider the scenario for one driver node, we follow this setting. The network size is chosen as 10, 20, 40, 60, 80, 100 accordingly. For all cases, we can see that the method we employed generates much more precise $\underline{E}$ compared to the existed tools.}\label{fig2}
\end{figure}

In the literature, it is common to use the trace of $\G$ to estimate the maximum eigenvalue of $\G$ \cite{Yan2012PRL,Li2017ConEng}. For the lower bound of $E$, we make a comparison of the precision between the existing result and the result obtained in this paper. From Fig.~\ref{fig2}, we find that the lower bounds derived in this paper are more exact.

By (\ref{EB}) with (\ref{upM}) and (\ref{lowM}), we have
\begin{equation}\label{upE}
\overline{E}\approx f(\underline{\alpha}, \underline{\beta}),
\end{equation}
and
\begin{equation}\label{lowE}
\underline{E}\approx \frac{1}{f(\overline{\alpha}, \overline{\beta})}.
\end{equation}
With only one driver node, we denote the node $h$ as the sole driver node with $b_{h1}=1$ and $b_{i1}=0 (i\neq h).$ Since $m_{ij}=q_{ij} f_{ij}$ and  $q_{ij}=p_{hi}p_{hj}$, we obtain
$
m_{ij}=\frac{p_{hi}p_{hj}}{\lambda_i+\lambda_j}(\text{e}^{(\lambda_i+\lambda_j)t_f}-1).
$
Furthermore, we have
$
\M^2(i,i)=\sum^n_{k=1}\frac{p^2_{hk}p^2_{hi}}{(\lambda_k+\lambda_i)^2}(\text{e}^{(\lambda_k+\lambda_i)t_f}-1)^2
$
and
$
\M^4(i,i)=\sum^n_{l=1}\left[\sum^n_{k=1}\frac{p^2_{hk}p_{hi}p_{hl}}{(\lambda_k+\lambda_i)(\lambda_k+\lambda_l)} \right.$ $\left.(\text{e}^{(\lambda_k+\lambda_i)t_f}-1)(\text{e}^{(\lambda_k+\lambda_l)t_f}-1)\right]^2.
$
Note that $\text{trace}(\L^2)=\|\L\|_F$ for arbitrary square matrix $\L$.
Then, we get the values of $\overline{\alpha}$ and $\overline{\beta}$ as
\begin{equation}\label{overalpha}
\overline{\alpha}=\text{trace}(\M^2)=\sum^n_{i=1}\sum^n_{k=1}\frac{p^2_{hk}p^2_{hi}}{(\lambda_k+\lambda_i)^2}(\text{e}^{(\lambda_k+\lambda_i)t_f}-1)^2,
\end{equation}
and
\begin{equation}\label{overbeta}
\overline{\beta}=\text{trace}(\M^4)=\sum^n_{i=1}\sum^n_{l=1}\left[\sum^n_{k=1} \frac{p^2_{hk}p_{hi}p_{hl}}{(\lambda_k+\lambda_i)(\lambda_k+\lambda_l)}(\text{e}^{(\lambda_k+\lambda_i)t_f}-1) (\text{e}^{(\lambda_k+\lambda_l)t_f}-1)\right]^2.
\end{equation}

\begin{figure}
\centering
\includegraphics[width=1\textwidth]{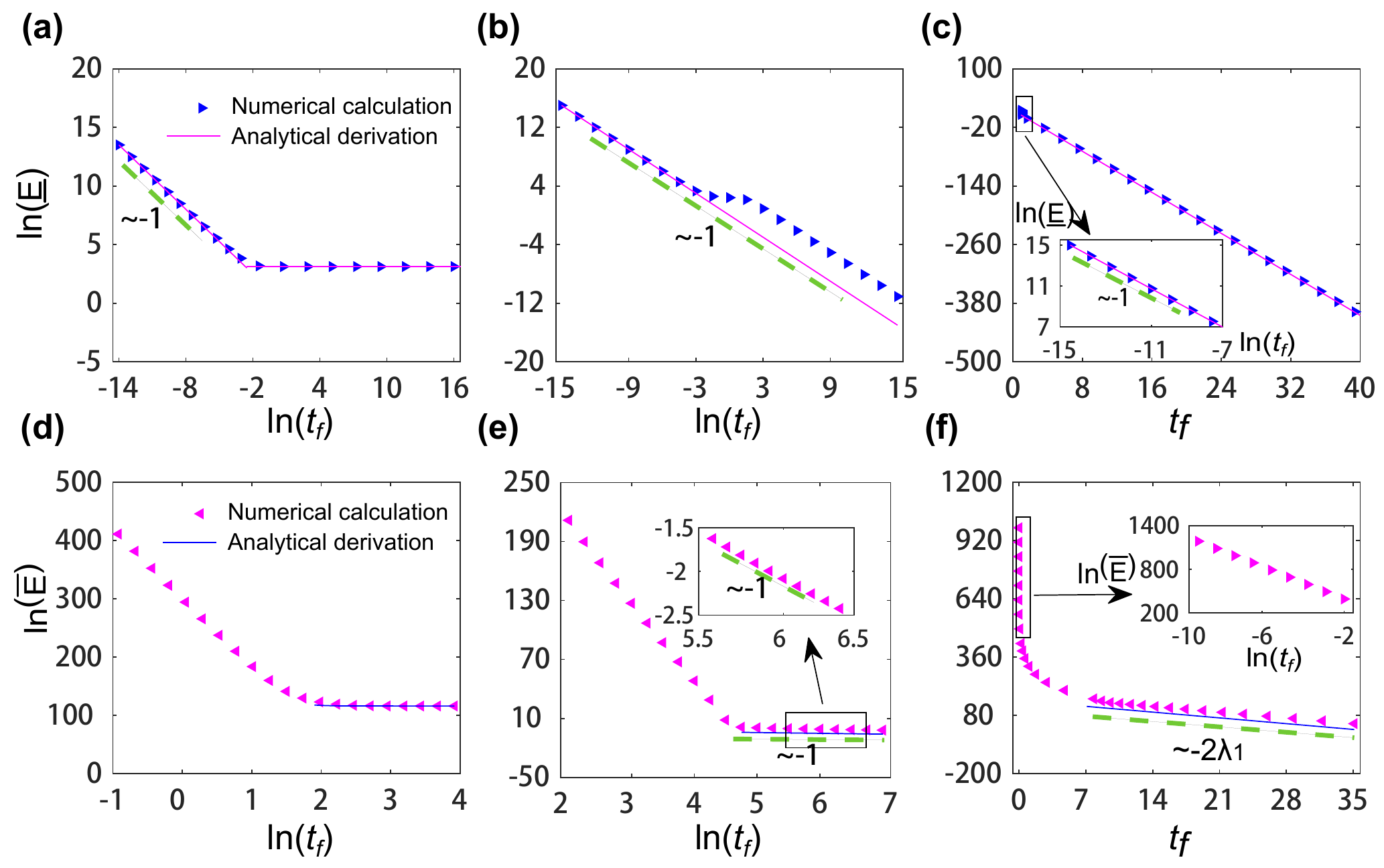}
\caption{The lower and upper bounds of energy for one driver node. The scaling behavior of the lower and upper bounds of energy cost is given for one driver node, and the summation of analytical results are presented in Tables \ref{tb5} and \ref{tb6}. In (a)-(c), with small $t_f$, $\underline{E}\sim t^{-1}_f$ for all $\A$. In (d)-(f) for upper bound, the slope of triangular trajectory is much less than $-1$. Parameters are selected the same as those given in Fig.~\ref{fig1}. The interval of the uniform distribution is $[0, 1]$ in (a)-(c), $[1, 3]$ in (d),  $[-1,0]$ in (e), and $[-5, -2]$ in (f). In (a), $a=-5$, by which $\A$ is ND with eigenvalues in $[-14.0266, -5]$. Similarly, in (b) and (e), $a=0$ such that $\A$ is NSD and PSD, respectively. In (c) and (d), $a=5$ such that $\A$ is ID. In (f), $a=3$, such that the minimum eigenvalue of $\A$ is $3$.
}
\label{fig4}
\end{figure}

Based on Eqs.~(\ref{overalpha}) and (\ref{overbeta}), we have discussed  and calculated the parameters $\overline{\alpha}$ and $\overline{\beta}$ in different cases (see Supplementary Information Sec.~S3). Accordingly, the upper and lower bounds of energy cost are given in Tables S1 and S2 in SI, and numerical validations of our analytical results are shown in Fig.~\ref{fig4}.

\subsection{$d$ driver nodes}\label{nd=d}
\begin{figure*}[!http]
\centering
\includegraphics[width=1\textwidth]{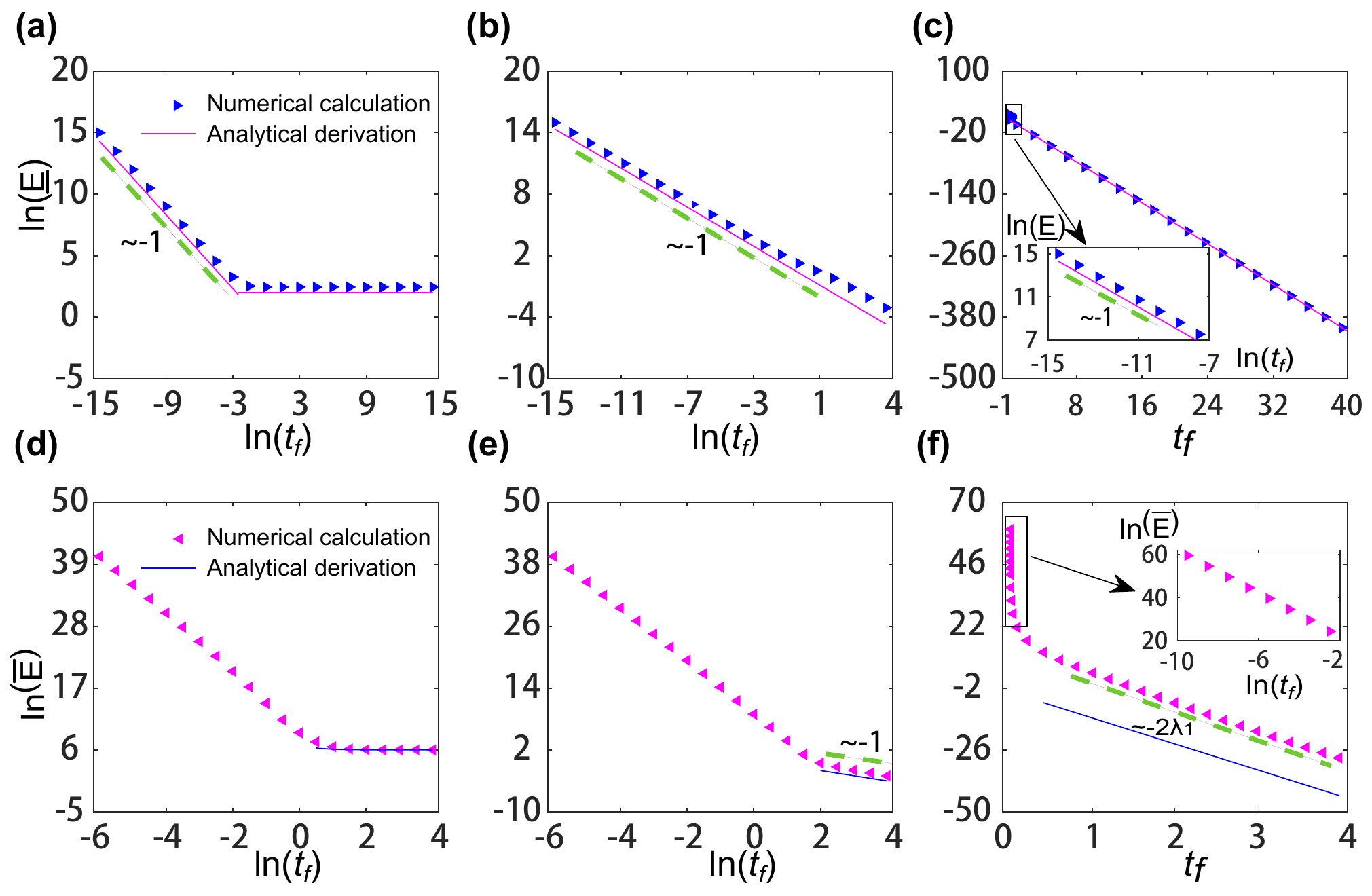}
\caption{The lower and upper bounds of control energy for $20$ driver nodes. In (a)-(c), with small $t_f$, $\underline{E}\sim t^{-1}_f$ for all $\A$. In (d)-(f) for upper bound, the slope of triangular trajectory is much less than $-1$. The summation of the analytical results are presented in Tables \ref{tb5} and \ref{tb6}. Parameters are selected as those given in Fig.~\ref{fig1}. The interval of uniform distribution is $[0, 1]$ in (a)-(d), and $[-1,0]$ in (e)-(f). In (a), $a=-5$, by which $\A$ is ND with eigenvalues in $[-12.5048, -5]$. Similarly, in (b) and (e), $a=0$ such that $\A$ is NSD and PSD, respectively. In (c) and (d), $a=5$ such that $\A$ is ID. Similarly, $a=5$ such that $\A$ is PD.
}\label{fig5}
\end{figure*}

In the case of $d$ driver nodes, we label them $m_1, m_2, \dots, m_d$. Hence $\B=[e_{m_1}, e_{m_2}, \dots, e_{m_d}]\in R^{n\times d}$, where $e_i=(0\,\dots\, 0\,\, 1\,\, 0\,\, \dots\,\, 0)^\text{T}\in R^n$ with all elements as $0$, except $i$th element as $1$. Let $\p_1=\B^\text{T}\p$, where $\p_1$ is a $d\times n$ matrix constituted by the rows $m_1$, $m_2$, $\dots$, $m_d$ of $\p$. Thus $\Q=\p^\text{T}_1\p_1$ with $q_{ij}=\sum^d_{k=1}p_{m_k i}p_{m_k j}$. By comparing the form of $m_{ij}=q_{ij}f_{ij}$ between the cases of one driver node and $d$ driver nodes, we find that only the form of $q_{ij}$ is different.
Therefore, in subsequent analysis and calculation, we can refer to the  Sec.~\ref{nd=1} to derive $\overline{\alpha}$ and $\overline{\beta}$ (see Sec.~S4 in SI for details). We summarize the lower bound of energy under $d$ driver nodes for different scenarios in Table S3 and the corresponding numerical validations are presented in Fig.~\ref{fig5}. In addition, the upper bound of energy is presented in Table S4.

\section{Discussion}
In this paper, we have investigated the scaling behavior of the bounds of minimum control energy for controlling complex networks in terms of the time given to achieve control. The bounds of minimum energy is determined by the maximum and the minimum eigenvalues of $\G$. The maximum eigenvalue is usually approximated by the trace of $\G$, while the approximation of the minimum eigenvalue has not yet been discussed in the existing literature. Here, we employ an effective method which not only provides more precise analytical expression than the trace for the approximation of the maximum eigenvalue, but also tells the analytical form of the minimum eigenvalues. All the derived theoretical laws are confirmed by numerical simulations.

Our framework also applies to weighted directed networks. When system (\ref{sys1}) is controllable, the matrix $\G$ is positive definite. When $\A$ is asymmetrical for directed networks, we can still obtain the specific form of $\G$. Based on $\G$, the lower bound of energy cost can be calculated by Eq.~(\ref{lowE}) with the traces of $\G^2$ and $\G^4$. For the upper bound of energy cost, we can apply the method to get the scaling behavior of energy by solving the inverse of $\G$ (see Sec.~S3 in SI).

Although natural systems are believed to operate with nonlinear dynamics, the type of nonlinearity and empirical  parameterization are usually hard to detect, especially for large systems. Besides, the generality of results cannot be guaranteed for some specific nonlinear systems. In contrast, the linear dynamics we analyzed here allows us to derive the theoretical insights, which is suitable for analyzing various complex networks. Even that we only consider static complex networks, our framework can also be employed to derive bounds of energy cost for controlling temporal networks by virtue of the effective matrix given in \cite{Li2017}. Specifically, utilizing estimations of the maximum and the minimum eigenvalues and some approximation techniques introduced in this paper, the \text{scaling} of energy for controlling temporal networks can be obtained.

\section*{Acknowledgements}
This work is supported by the National Natural Science Foundation of China (NSFC) under grants no. 61751301 and no. 61533001.
A.L. acknowledgements the
Human Frontier Science Program Postdoctoral Fellowship (Grant: LT000696/2018-C),
the generous support from Foster Lab at Oxford, and the Chair of Systems Design at ETH Z\"urich. ETH Z\"urich, Weinbergstrasse 56/58, Z\"urich CH-8092, Switzerland.
G.Z. acknowledgements the financial support from the Hong Kong Research Grant council (RGC) grants (No. 15206915, No. 15208418).
%\begin{enumerate}
%\item  distribution has been given, but only under infinite time which does not hold for empirical scenario.
%\item directed networks
%\item  unstable systems
%\item  temporal networks and nonlinear dynamics
%\end{enumerate}

%\section{To do list}
%
%\begin{enumerate}
%\item fix the framework to show results (according to the main findings)
%\item draw figures (first with cartoon, later for each specific scenario)
%\item show numerical calculations for large systems (empirical networks)
%\end{enumerate}

%\bibliography{../Public/bibliography}

\begin{table}[!http]
\caption{The lower bound of control energy $\underline{E}$. No matter how many driver nodes there are, for small $t_f$, $\underline{E}\sim t^{-1}_f$. For large $t_f$, when $\A$ is ND (negative definite), $\underline{E}$ approaches to a constant irrespective of $t_f$, ($C_1$ for one driver node, $C_2$ for $d$ driver nodes and $2|\lambda_n|$ for $n$ driver nodes), where $C_1$ and $C_2$ are given as Eq.~(\ref{lowE}) with Eqs.~(S6) (S7) in Sec.~S3 and with Eqs.~(S45) (S46) in Sec.~S4 of SI, respectively. When $\A$ is NSD (negative semi-definite) with large $t_f$, $\underline{E}\approx t^{-1}_f$ under $1$ and $n$ driver nodes; while it approaches $t^{-1}_f$(detailed forms are given as Eq.~(\ref{lowE}) with Eqs.~(S47) and (S48) in SI). In addition, when $\A$ is not ND (including the cases of indefinite, positive semi-definite, and positive definite), $\underline{E}\sim \text{e}^{-2\lambda_nt_f}$ holds for large $t_f$.}\label{tb5}
\centering
\fontsize{8}{15}\selectfont
\begin{tabular}{cc|c|c|c}
\hline
\multicolumn{2}{c|}{Number of driver nodes}&$1$&$d$&$n$\\
\hline
\multicolumn{2}{c|}{Small $t_f$}&$ t_f^{-1}$&$\sim t_f^{-1}$&$t^{-1}_f$\\
%\hline
\multirow{3}{*}{Large $t_f$}&ND&$C_1$&$C_2$&$2|\lambda_n|$\\
%\hline
&NSD&$ t_f^{-1}$& $\sim t_f^{-1}$&$t^{-1}_f$\\
%\hline
&Not ND&$\sim \text{e}^{-2\lambda_n t_f}$&$\sim \text{e}^{-2\lambda_n t_f}$&$\sim \text{e}^{-2\lambda_nt_f}$\\
\hline
\end{tabular}
\label{table1}
\end{table}

\begin{table}[!http]
\caption{
The upper bound of control energy $\overline{E}$. For small $t_f$, both $N_0-N_{\min}$ and $N'_0-N'_{\min}$ are much larger than $1$, where the detailed meanings of $N_0, N_{\min}, N'_0$ and $N'_{\min}$ are given in Secs.~S3 and S4 of SI. For large $t_f$, when $\A$ is PD (positive definite), $\overline{E}\sim \text{e}^{-2\lambda_1t_f}$ for arbitrary number of driver nodes; when $\A$ is PSD (positive semi-definite), $\overline{E}\sim t^{-1}_f$; when $\A$ is not PD (including the cases of indefinite, negative semi-definite, and negative definite), $\overline{E}$ approaches to a constant irrespective of the magnitude of $t_f$ for large $t_f$ ($C_3$ for one driver node, $C_4$ for $d$ driver nodes, and $2|\lambda_1|$ for $n$ driver nodes), where $C_3$ has different forms for different $\A$ (detailed forms are presented in Table~S2 of Sec.~S3 of SI). }\label{tb6}
\centering
\fontsize{8}{15}\selectfont
\begin{tabular}{cc|c|c|c}
\hline
\multicolumn{2}{c|}{Number of driver nodes}&$1$&$d$&$n$\\
\hline
\multicolumn{2}{c|}{Small $t_f$}&$\sim t_f^{-(N_0-N_{\min})/2}$&$\sim t_f^{-(N'_0-N'_{\min})/2}$&$t^{-1}_f$\\
%\hline
\multirow{3}{*}{Large $t_f$}&PD&$\sim \text{e}^{-2\lambda_1t_f}$&$\sim \text{e}^{-2\lambda_1t_f}$&$\sim \text{e}^{-2\lambda_1t_f}$\\
%\hline
&PSD&$\sim t^{-1}_f$&$\sim t^{-1}_f$&$t^{-1}_f$\\
%\hline
&Not PD&$C_3$&$C_4$&$2|\lambda_1|$\\
\hline
\end{tabular}
\end{table}


\end{document}